\documentstyle[twoside,amstex]{article}
\input amssym.def
\input amssym.tex
\def\bc{\begin{center}}
\def\ec{\end{center}}
\def\no{\noindent}

\begin{document}
\thispagestyle{empty} \vspace*{3 true cm} \pagestyle{myheadings}
\markboth {\hfill {\sl Huanyin Chen and Marjan Sheibani}\hfill}
{\hfill{\sl Periodicity and J-Clean-like Rings
}\hfill} \vspace*{-1.5 true cm} \bc{\large\bf Periodicity and J-Clean-like Rings}\ec

\vskip10mm\bc{{\bf Huanyin Chen}\\[2mm]
Department of Mathematics, Hangzhou Normal University\\
Hangzhou 310036, China\\
huanyinchen@@aliyun.com}\ec

\vskip2mm \bc{{\bf Marjan
Sheibani}\\[2mm]
Faculty of Mathematics, Statistics and Computer Science\\
Semnan University, Semnan, Iran\\
m.sheibani1@@gmail.com}\ec

\vskip10mm
\begin{minipage}{120mm}
\no {\bf Abstract:} A ring $R$ is periodic provided that for any
$a\in R$ there exist distinct elements $m,n\in {\Bbb N}$ such that
$a^m=a^n$. We shall prove that periodicity is inherited by all generalized matrix rings.
A ring $R$ is called strongly periodic if for any $a\in R$ there
exists a potent $p\in R$ such that $a-p$ is in its prime radical and $ap=pa$.
A ring $R$ is J-clean-like if for
any $a\in R$ there exists a potent $p\in R$ such that $a-p$ is in its Jacobson radical.
Furthermore, we completely determine the connections between strongly periodic rings and
periodic rings. The relations among J-clean-like rings and these rings are also obtained.


\end{minipage}

\vskip10mm\section{Introduction} \vskip4mm A ring $R$ is periodic
provided that for any $a\in R$ there exist distinct elements
$m,n\in {\Bbb N}$ such that $a^m=a^n$. Examples of periodic rings
are finite rings and Boolean rings.
There are many interesting problems related to periodic rings. We
explore, in this article, the periodicity of a type of generalized
matrix rings. An element $p\in R$ is potent if $p=p^m$ for some
$m\geq 2$. For later convenience we state here some elementary
characterizations of periodic rings:

\vskip4mm \hspace{-1.8em} {\bf Theorem 1.1.}\ \ {\it Let $R$ be a
ring. Then the following are equivalent:}
\begin{enumerate}
\item [(1)]{\it $R$ is periodic.}
\vspace{-.5mm}
\item [(2)]{\it For any $a\in R$, there exists some $m\geq 2$ such that $a^m=a^{m+1}f(a)$ for some $f(t)\in {\Bbb Z}[t]$.}
\item [(3)]{\it For any $a\in R$, there exists some $m\geq 2$ such that $a-a^m\in R$ is nilpotent.}
\vspace{-.5mm}
\item [(4)]{\it For any $a\in R$, there exists a potent $p\in R$ such that $a-p\in R$ is nilpotent and $ap=pa$.}
\end{enumerate}
Here, the equivalences of all items are stated in ~\cite[Lemma 2]{A}, \cite[Proposition 2]{Cha}, \cite[Theorem 3]{MO}, and the simple implication from $(3)$ to $(2)$.
A Morita context $(A,B,M,N,\psi,\varphi)$ consists of two rings
$A$ and $B$, two bimodules $_AN_B$ and $_BM_A$, and a pair of
bimodule homomorphisms $\psi : N\bigotimes\limits_{B}M\to A$ and
$\varphi: M\bigotimes\limits_{A} N\to B$ which satisfy the
following associativity: $\psi \big(n\bigotimes m\big)n'=n\varphi
\big(m\bigotimes n'\big)$ and $\varphi \big(m\bigotimes
n\big)m'=m\psi \big(n\bigotimes m'\big)$ for any $m,m'\in M,
n,n'\in N$. These conditions ensure that the set $T$ of
generalized matrices $\left(
\begin{array}{cc}
a&n\\
m&b
\end{array}
\right); a\in A,b\in B,m\in M,n\in N$ will form a ring with addition defined componentwise and
with multiplication defined by
$$\left(
\begin{array}{cc}
a_1&n_1\\
m_1&b_1
\end{array}
\right)\left(
\begin{array}{cc}
a_2&n_2\\
m_2&b_2
\end{array}
\right)=\left(
\begin{array}{cc}
a_1a_2+\psi(n_1\bigotimes m_2)&a_1n_2+n_1b_2\\
m_1a_2+b_1m_2&\varphi(m_1\bigotimes n_2)+b_1b_2
\end{array}
\right),$$ called the
ring of the Morita context(cf. ~\cite{TZ}). The class of rings of the Morita
contexts is a type of generalized matrix rings. For instances, all
$2\times 2$ matrix rings and all triangular matrix rings.

Let $T$ be the ring of a Morita context $(A,B,$ $M,N,$
$\psi,\varphi)$. We prove, in Section 2, that if $im(\psi)$ and
$im(\varphi)$ are nilpotent, then $A$ and $B$ are periodic if and
only if so is $T$. This provides a large new class of periodic
rings for generalized matrix rings.

It is an attractive problem to express an element in a ring as the
sum of idempotents and units (cf.~\cite{B}, \cite{CH}, \cite{D} and
\cite{H}). We say that a ring $R$ is clean provided that every
element in $R$ is the sum of an idempotent and a unit. Such rings
have been extensively studied recently years, see ~\cite{CH1} and
~\cite{T}. This motivates us to combine periodic rings with clean
rings together, and investigate further properties of related
rings.

For a ring $R$ the prime radical is denoted by $P(R)$, i.e., $P(R)$ is the intersection of all prime ideals of $R$. We now introduce a new type of rings. A ring $R$
is said to be strongly periodic provided that for any $a\in R$ there
exists a potent $p\in R$ such that $a-p\in P(R)$ and $ap=pa$.
Strongly periodic rings form a subclass of periodic rings. We
shall prove that a ring $R$ is strongly periodic if and only if
for any $a\in R$ there exists a potent $p\in R$ such that $a-p\in
P(R)$, and determine completely the connections between these ones
and periodic rings. A ring is 2-primal provided that its
prime radical coincides with the set of nilpotent elements of the
ring. It is proved that a ring $R$ is strongly periodic if and only
if $R$ is a 2-primal periodic ring. From this, we show that the
strong periodicity will be inherited by generalized matrix rings.

Replacing the prime radical $P(R)$ by the Jacobson radical $J(R)$,
we introduce a type of rings which behave like that of periodic
rings. We say that a ring $R$ is J-clean-like provided that for
any $a\in R$ there exists a potent $p\in R$ such that $a-p\in
J(R)$. This is a natural generalization of J-clean rings
~\cite{CH}. Many properties of periodic rings are extended to
these ones. We shall characterize J-clean-like rings and obtain
the relations among these rings.

Throughout, all rings are associative with an identity. $M_n(R)$
will denote the ring of all $n\times n$ matrices over $R$
with an identity $I_n$. $N(R)$ stands for the set of all nilpotent
elements in $R$. $C(R)$ denote the center of $R$. $P(R)$ and
$J(R)$ denote the prime radical and Jacobson radical of $R$,
respectively.

\section{Periodic rings}

\vskip4mm The purpose of this section is to investigate the
periodicity for Morita contexts. The following lemma is
known~\cite[Lemma 3.1.23]{MB}, and we include a simple proof for
the sake of completeness.

\vskip4mm \hspace{-1.8em} {\bf Lemma 2.1.}\ \ {\it A ring $R$ is
periodic if and only if for any $a,b\in R$, there exists an $n\in
{\Bbb N}$ such that $a-a^n,b-b^n\in N(R)$.}
\vskip2mm\hspace{-1.8em} {\it Proof.}\ \ $\Longleftarrow$ For any $a\in R$, we can find $n\in {\Bbb N}$ such that $a-a^n\in N(R)$. This implies that $R$ is periodic, by Theorem 1.1.

$\Longrightarrow $ Suppose that $R$ is periodic. For any
$a,b\in R$, we can find $p,q,s,t\in R$ $(p<q, s<t)$ such that
$a^p=a^q$ and $b^s=b^t$. Hence, $a^{ps}=a^{qs}$ and
$b^{ps}=b^{pt}.$ This implies that
$$a^{ps}=a^{ps}a^{(q-p)s}=a^{ps}a^{2(q-p)s}=\cdots =a^{ps}a^{(t-s)p(q-p)s}.$$
Likewise, we get $b^{ps}=b^{ps}b^{(q-p)s(t-s)p}.$ Choose $k=ps$
and $l=ps+(t-s)p(q-p)s$. Then $a^k=a^l, b^k=b^l$ $(k<l)$.
Thus, $a^k=a^l=a^{(l-k)+k}=\cdots
=a^{k(l-k)+k}$, and so $a^k=\big(a^k\big)^{l-k+1}$. This implies
that
$\big(a^{k(l-k)}\big)^2=\big(a^{k(l-k)+k}\big)\big(a^{k(l-k)-k}\big)=a^k\big(a^{k(l-k)-k}\big)=a^{k(l-k)}.$
Choose $n=k(l-k)$. Then
$\big(a-a^{n+1}\big)^n=a^n\big(1-a^n\big)^n=a^n\big(1-a^n\big)=0$.
Thus, $a-a^n\in N(R)$. Likewise, $b-b^n\in N(R)$. Therefore we complete the proof.\hfill$\Box$

\vskip4mm \hspace{-1.8em} {\bf Theorem 2.2.}\ \ {\it Let $T$ be
the ring of a Morita context $(A,B,$ $M,N,$ $\psi,\varphi)$. If
$im(\psi)$ and $im(\varphi)$ are nilpotent, then $A$ and $B$ are
periodic if and only if so is $T$.} \vskip2mm\hspace{-1.8em} {\it
Proof.}\ \ Suppose $A$ and $B$ are periodic. For any $\left(
\begin{array}{cc}
a&n\\
m&b
\end{array}
\right)\in T$, as in the proof of Lemma 2.1, there exists a $k\in
{\Bbb N}$ such that $a-a^k\in N(A)$ and $b-b^k\in N(B)$. Hence,
$$\left(
\begin{array}{cc}
a&n\\
m&b
\end{array}
\right)-\left(
\begin{array}{cc}
a&n\\
m&b
\end{array}
\right)^k=\left(
\begin{array}{cc}
a-a^k+c&*\\
**&b-b^k+d
\end{array}
\right),$$ where $c\in im(\psi)$ and $d\in im(\varphi)$. Write
$(a-a^k)^l=0$ and $(b-b^k)^l=0$. By hypothesis, $im(\psi)$ and
$im(\varphi)$ are nilpotent ideals of $A$ and $B$, respectively.
Say $\big(im(\psi)\big)^s=0$ and $\big(im(\varphi)\big)^t=0$.
Choose $p=\max (s,t)$ and $q=p(l+1)$. Then
$$\big(a-a^k+c\big)^q=0~\mbox{and}~\big(b-b^k+d\big)^q=0.$$
Obviously, $$\left(
\begin{array}{cc}
a-a^k+c&*\\
**&b-b^k+d
\end{array}
\right)^{q+1}\in \left(
\begin{array}{cc}
im(\psi)&N\\
M&im(\varphi)
\end{array}
\right).$$
Set $NM:=im(\psi)$ and $MN:=im(\varphi)$. We see that
$$\left(
\begin{array}{cc}
NM&N\\
M&MN
\end{array}
\right)^2\subseteq \left(
\begin{array}{cc}
NM&(NM)N\\
(MN)M&MN
\end{array}
\right).$$ For any $l\in {\Bbb N}$, by induction, one easily checks that
$$\left(
\begin{array}{cc}
NM&N\\
M&MN
\end{array}
\right)^{2l}\subseteq \left(
\begin{array}{cc}
NM&(NM)N\\
(MN)M&MN
\end{array}
\right)^l\subseteq \left(
\begin{array}{cc}
(NM)^l&(NM)^lN\\
(MN)^lM&(MN)^l
\end{array}
\right).$$
Choose $j=2p(q+1)$. As $(NM)^p=(MN)^p=0$, we get
$$\left(
\begin{array}{cc}
a-a^k+c&*\\
**&b-b^k+d
\end{array}
\right)^{j}=\left(
\begin{array}{cc}
0&t\\
s&0
\end{array}
\right)$$ for some $s\in M,t\in N$.
Hence, $$\left(
\begin{array}{cc}
0&t\\
s&0
\end{array}
\right)^2=\left(
\begin{array}{cc}
\psi(t\bigotimes s)&0\\
0&\varphi(s\bigotimes t)
\end{array}
\right),$$ and so
$$\big(\left(
\begin{array}{cc}
a&n\\
m&b
\end{array}
\right)-\left(
\begin{array}{cc}
a&n\\
m&b
\end{array}
\right)^k\big)^{2jp}=0.$$ Accordingly, $T$ is periodic, by Theorem
1.1. The converse is obvious.\hfill$\Box$

\vskip4mm Let $R$ be a ring, and let $s\in C(R)$. Let
$M_{(s)}(R)=\{ \left(
\begin{array}{cc}
a&b\\
c&d
\end{array}
\right)~|~a,b,c,d\in R\}$, where the operations are defined as
follows:
$$
\begin{array}{c}
\left(
\begin{array}{cc}
a&b\\
c&d
\end{array}
\right)+\left(
\begin{array}{cc}
a'&b'\\
c'&d'
\end{array}
\right)=\left(
\begin{array}{cc}
a+a'&b+b'\\
c+c'&d+d'
\end{array}
\right),\\
\left(
\begin{array}{cc}
a&b\\
c&d
\end{array}
\right)\left(
\begin{array}{cc}
a'&b'\\
c'&d'
\end{array}
\right)=\left(
\begin{array}{cc}
aa'+sbc'&ab'+bd'\\
ca'+dc'&scb'+dd'
\end{array}
\right). \end{array}$$ Then $M_{(s)}(R)$ is a ring with the
identity $\left(
\begin{array}{cc}
1_R&0\\
0&1_R
\end{array}
\right)$. Recently, the strong cleanness of such type generalized
matrix rings was studied in ~\cite{T}. For the periodicity of such
rings, we derive

\vskip4mm \hspace{-1.8em} {\bf Corollary 2.3.}\ \ {\it Let $R$ be
periodic, and let $s\in N(R)\bigcap C(R)$. Then $M_{(s)}(R)$ is
periodic.} \vskip2mm\hspace{-1.8em} {\it Proof.}\ \ Let $\psi:
R\otimes R\to R, n\otimes m\mapsto snm$ and $\varphi: R\otimes
R\to R, m\otimes n\mapsto smn$. Then
$M_s(R)=(R,R,R,R,\psi,\varphi)$. As $s\in N(R)\bigcap C(R)$, we
see that $im(\varphi)$ and $im(\psi)\subseteq J(R)$ are nilpotent,
and we are through by Theorem 2.2.\hfill$\Box$

\vskip4mm As a consequence, a ring $R$ is periodic if and only if
so is the trivial Morita context $M_{(0)}(R)$. Choosing $s=0\in R$, we are through from
Corollary 2.3. Given a ring $R$ and an $R$-$R$-bimodule $M$, the
trivial extension of $R$ by $M$ is the ring $T(R,M)= R\oplus M$
with the usual addition and the following multiplication: $(r_1,
m_1)(r_2, m_2) = (r_1r_2, r_1m_2 + m_1r_2)$.

\vskip4mm \hspace{-1.8em} {\bf Corollary 2.4.}\ \ {\it Let $R$ be
a ring, and let $M$ be a $R$-$R$-bimodule. Then the following are
equivalent:}
\begin{enumerate}
\item [(1)]{\it $R$ is periodic.}
\vspace{-.5mm}
\item [(2)]{\it $T(R,M)$ is periodic.}
\end{enumerate}
\vspace{-.5mm} {\it Proof.}\ \  $(1)\Rightarrow (2)$ Let $R$ be a
periodic ring and let $S=\left(
\begin{array}{cc}
R&M\\
0&R
\end{array}
\right)$. It is obvious by Theorem 2.2 that $S$ is periodic.
Clearly, $T(R,M)$ is a subring of $S$, and so proving $(2)$.

$(2)\Rightarrow (1)$ Let $T(R,M)$ be a periodic ring. As $R$ is
isomorphic to a subring of $T(R,M)$, and so $R$ is periodic.
\hfill$\Box$

\vskip4mm \hspace{-1.8em} {\bf Example 2.5.}\ \ Let $R$ be
periodic, let $$A=B=\left(
\begin{array}{ccc}
R&0&0\\
0&R&0\\
0&0&R
\end{array}
\right),M=\left(
\begin{array}{ccc}
0&0&0\\
0&0&0\\
0&R&0
\end{array}
\right)~\mbox{and}~N=\left(
\begin{array}{ccc}
0&0&0\\
0&0&0\\
R&R&0
\end{array}
\right),$$ and let $\psi : N\bigotimes\limits_{B}M\to A,
\psi(n\otimes m)=nm$ and $\phi :M\bigotimes\limits_{A} N\to B,
\phi(m\otimes n)=mn$. Then $T=(A,B,M,N,\psi ,\phi )$ is a Morita
context with zero pairings, i.e., $T$ is a trivial Morita context. Hence, $im(\psi)$ and $im(\varphi)$
are nilpotent. Clearly, $A$ and $B$ are both periodic. In light of
Theorem 2.2, $T$ is periodic.

\vskip4mm Let $R$ be a ring, and let $\alpha$ be an endomorphism
of $R$. Let $T_n(R,\alpha)$ be the set of all upper triangular
matrices over the rings $R$. For any $(a_{ij}),(b_{ij})\in
T_n(R,\alpha)$, we define $(a_{ij})+(b_{ij})=(a_{ij}+b_{ij})$, and
$(a_{ij})(b_{ij})=(c_{ij})$ where
$c_{ij}=\sum\limits_{k=i}^{n}a_{ik}\alpha^{k-i}\big(b_{kj}\big)$.
Then $T_n(R,\alpha)$ is a ring under the preceding addition and
multiplication (cf.~\cite{JM}). Clearly, $T_n(R,\alpha)$ will be $T_n(R)$ only
when $\alpha$ is the identity morphism.

\vskip4mm \hspace{-1.8em} {\bf Lemma 2.6.}\ \ {\it Let $R$ be
periodic, and let $\alpha: R\to R$ be an endomorphism. Then}
\begin{enumerate}
\item [(1)]{\it $R[[x,\alpha]]/(x^n)$ is periodic.}
\vspace{-.5mm}
\item [(2)]{\it $T_n(R,\alpha)$ is periodic for all $n\in {\Bbb
N}$.}
\end{enumerate}
\vspace{-.5mm} {\it Proof.}\ \ $(1)$ For any $f(x)\in
R[[x]]/(x^n)$, there exists an $m\in {\Bbb N}$ such that
$f(0)-f^m(0)\in N(R)$. Hence, $f(x)-f^m(x)\in
N\big(R[[x]]/(x^n)\big)$. According to Theorem 1.1, $R[[x]]/(x^n)$
is periodic.

$(2)$ For any $(a_{ij})\in T_n(R,\alpha)$, as in the proof of
Lemma 2.1, we can find an $m\in {\Bbb N}$ such that
$a_{ii}-a_{ii}^m\in N(R)$ for each $i$. Thus,
$(a_{ij})-(a_{ij})^m\in N\big(T_n(R,\alpha)\big)$, as
required.\hfill$\Box$

\vskip4mm We are now ready to prove:

\vskip4mm \hspace{-1.8em} {\bf Theorem 2.7.}\ \ {\it Let $R$ be
periodic. Then $M_{(x^m)}\big(R[[x]]/(x^n)\big)$ is periodic for
all $1\leq m\leq n$.} \vskip2mm\hspace{-1.8em} {\it Proof.}\ \
Choose $\alpha=1$. Then $R[[x]]/(x^n)$ is periodic, by Lemma 2.6.
Choose $s=x^m (1\leq m\leq n)$. Then $s\in
N\big(R[[x]]/(x^n)\big)\bigcap C\big(R[[x]]/(x^n)\big)$. Applying
Corollary 2.3 to $R[[x]]/(x^n)$, $M_{(x^m)}\big(R[[x]]/(x^n)\big)$
is periodic, as asserted.\hfill$\Box$

\vskip4mm \hspace{-1.8em} {\bf Corollary 2.8.}\ \ {\it Let $R$ be
a finite ring. Then $M_{(x^m)}\big(R[[x]]/(x^n)\big)$ is periodic
for all $1\leq m\leq n$.}
 \vskip2mm\hspace{-1.8em} {\it Proof.}\ \ Since every finite ring is periodic, we complete the proof by
Theorem 2.7.\hfill$\Box$

\section{Strongly periodic rings}

\vskip4mm An element $a\in R$ is strongly nilpotent if for any sequence $x_0=a, x_1, x_2,\cdots $ with $x_{i}\in x_{i-1}Rx_{i-1}$ for each $i\in {\Bbb N}$, there exists some $n$ such that $x_n=0$. As is well known, $P(R)$ is exactly the set of all strongly nilpotent elements in $R$ (cf.~\cite{L}). Hence, $P(R)$ is a nil ideal of $R$, i.e., every element in $P(R)$ is nilpotent.
A ring $R$ is potent if for any $a\in R$, there exists some $n\geq 2$ such that $a=a^n$. An ideal $I$ of a ring $R$ is locally nilpotent if,
every finitely generated subring of elements belonging to $I$ is nilpotent. Clearly, an ideal $I$ of a ring $R$ is locally nilpotent if and only if
$RxR$ is nilpotent for any $x\in I$. Recall that $J(R)$ consists of all $x\in R$ such that $1+RxR$ is included in the set of units of $R$. We now derive

\vskip4mm \hspace{-1.8em} {\bf Theorem 3.1.}\ \ {\it Let $R$ be a ring. Then the following are equivalent:}
\begin{enumerate}
\item [(1)] {\it $R$ is
strongly periodic.}
\vspace{-.5mm}
\item [(2)] {\it $R$ is periodic and $N(R)$ is a locally nilpotent ideal of $R$.}
\vspace{-.5mm}
\item [(3)] {\it $R/J(R)$ is potent, every potent lifts modulo $J(R)$ and $J(R)$ is locally nilpotent.}
\end{enumerate}
\vspace{-.5mm} {\it Proof.}\ \ $(1)\Rightarrow (2)$ Clearly, $R$ is periodic. Let $x\in N(R)$. Then we can find a potent $p\in R$ such that $w:=x-p\in P(R)$. Write $x^n=0$ for some $n\in {\Bbb N}$. Then $p^n=(x-w)^n\in P(R)$. This shows that $p\in R$ is nilpotent, and so $p=0$; hence, $x=w\in P(R)$. We infer that $N(R)=P(R)$ is an ideal of $R$.

For any $x\in P(R)$,
We claim that $RxR$ is nilpotent. If not, set
$\Delta=\{R,RxR,(RxR)^2,\cdots \}$ and $\Omega_P=\{ RyR~|~y\in
P\}$ for any ideal $P$. Let
$$\Omega=\{P\mid P~\mbox{is an ideal of}~R~\mbox{such
that}~\Omega_P\bigcap \Delta=\emptyset\}.$$ Clearly, $\Omega$ is a
nonempty inductive set. By using Zorn's Lemma, we have an ideal
$Q$ which is maximal in $\Omega$. If $Q$ is not prime, then we can
find some $a,b\in R$ such that $aRb\subseteq Q$, while $a,
b\not\in Q$. Thus, there exist $m,n\in \mathbb{N}$ such that
$(RxR)^m\subseteq Q+RaR$ and $(RxR)^{n}\subseteq Q+RbR$. This
implies that $(RxR)^{m+n}\in Q+RaRbR\subseteq Q$, a contradiction.
We infer that $Q\in Spec(R)$, and so $x\in Q$, a
contradiction.

Thus, $N(R)$ is locally nilpotent.

$(2)\Rightarrow (1)$ Let $x\in N(R)$. As $N(R)$ is locally nilpotent, $RxR$ is nilpotent.
Write $(RxR)^m=0 (m\in {\Bbb N})$. Let $P$ be an arbitrary prime ideal of $R$. Then
$(RxR)^m\subseteq P$; and so $RxR\subseteq P$. This implies that $x\in P$; hence,
$N(R)\subseteq P(R)$. The implication is true, by Theorem 1.1.

$(1)\Rightarrow (3)$ For any
$a\in R$ there exists some potent $p\in R$ such that $a-p\in
P(R)\subseteq J(R)$. Hence, $\overline{a}=\overline{p}$ in
$R/J(R)$. Therefore $R/J(R)$ is potent.

Let $x\in J(R)$. Then there exists a potent $p\in R$ such
that $x-p\in P(R)$; hence, $p=x-(x-p)\in J(R)$. Write $p=p^{m}
(m\geq 2)$. then $p(1-p^{m-1})=0$, and so $p=0$. Hence, $x\in
P(R)$. By the preceding discussion, $RxR$ is nilpotent, and therefore $J(R)$ is
locally nilpotent.

$(3)\Rightarrow (1)$ Let $a\in R$. Then $a-a^n\in J(R)$ for some
$n\geq 2$. As $J(R)$ is locally nilpotent, it is nilpotent, and so $a-a^n\in N(R)$. In
view of Theorem 1.1, $R$ is periodic. Let $x\in N(R)$. Then
$\overline{x}\in R/J(R)$ is potent; hence,
$\overline{x}=\overline{0}$ in $R/J(R)$. That is, $x\in
J(R)$. By hypothesis, $J(R)$ is locally nilpotent; hence, $RxR$ is nilpotent. As in the proof in $(2)\Rightarrow (1)$, we see that
$x\in P(R)$. Thus, $N(R)\subseteq P(R)$.

For any $a\in R$, there exists a potent $p\in R$ and a $w\in N(R)$
such that $a=p+w$ and $pw=wp$, by Theorem 1.1. By the preceding
discussion, $w\in P(R)$. This proving $(1)$.\hfill$\Box$

\vskip4mm \hspace{-1.8em} {\bf Corollary 3.2.}\ \ {\it A ring $R$
is strongly periodic if and only if for any $a\in R$ there exists
a potent $p\in R$ such that $a-p\in P(R)$.} \vskip2mm\hspace{-1.8em}
{\it Proof.}\ \ $\Longrightarrow$ This is trivial.

$\Longleftarrow$ For any $a\in R$ there exists potent $p\in R$
such that $a-p\in P(R)\subseteq J(R)$. Hence, $R/J(R)$ is potent.
For any $x\in J(R)$, there exists a potent $q\in R$ such that
$x-q\in P(R)$. Hence, $q=x-(x-q)\in J(R)$. Write $q=q^m (m\geq
2)$. Then $q(1-q^{m-1})=0$, and so $q=0$. We infer that $x\in
P(R)$. As in the proof in Theorem 3.1, $RxR$ is nilpotent, and so $J(R)$ is locally
nilpotent. This result follows, by using Theorem 3.1.\hfill$\Box$

\vskip4mm Birkenmeier-Heatherly-Lee introduced the
concept of 2-primal ring $R$, i.e., $N(R)=P(R)$ (~\cite{Bi}). Shin proved that a ring is 2-primal if
and only if each of its minimal prime ideal is completely prime.
A ring $R$ is weakly periodic provided that for any $a\in R$ there
exists a potent $p\in R$ such that $a-p\in N(R)$ ~\cite{MB}. We now derive

\vskip4mm \hspace{-1.8em} {\bf Theorem 3.3.}\ \ {\it A ring $R$ is
strongly periodic if and only if $R$ is a 2-primal weakly periodic
ring.} \vskip2mm\hspace{-1.8em} {\it Proof.}\ \ $\Longrightarrow $
Clearly, $R$ is weakly periodic. For any $a\in N(R)$, there exists
a potent $p\in R$ such that $w:=a-p\in P(R)$. Hence, $p=a-w$.
Write $a^{m}=0 (m\in {\Bbb N})$. Then $p^m\in P(R)$, and so $p\in
N(R)$. This implies that $p=0$, and so $a=w\in P(R)$. Thus,
$N(R)=P(R)$, and so $R$ is 2-primal.

$\Longleftarrow $ Let $a\in R$. Since $R$ is weakly periodic,
there exists a potent $p\in R$ such that $a-p\in N(R)$. As $R$ is
2-primal, $N(R)\subseteq P(R)$, we get $a-p\in P(R)$. Therefore we
complete the proof, by Corollary 3.2.\hfill$\Box$

\vskip4mm A ring $R$ is called strongly 2-primal provided that
$R/I$ is 2-primal for all ideals $I$ of $R$.

\vskip4mm \hspace{-1.8em} {\bf Corollary 3.4.}\ \ {\it A ring $R$
is strongly periodic if and only if the following two conditions hold:}
\begin{enumerate}
\item [(1)] {\it $R$ is weakly periodic;}
\vspace{-.5mm} \vspace{-.5mm}
\item [(2)] {\it Every prime ideal of $R$ is completely prime.}
\end{enumerate}
\vspace{-.5mm} {\it Proof.}\ \ $\Longrightarrow$ $(1)$ is obvious.
Clearly, $R/P(R)$ is potent. As in well known, every potent ring
is commutative (cf.~\cite[Theorem 1 in Chapter X]{J}), and so $R/P(R)$ is commutative. Hence, $R/P(R)$ is
strongly 2-primal. In view of ~\cite[Proposition 1.2]{Lee}, every
prime ideal of $R$ is completely prime.

$\Longleftarrow$ In view of ~\cite[Proposition 1.2]{Lee}, $R$ is
strongly 2-primal, and then it is 2-primal. This completes the
proof, in terms of Theorem 3.3.\hfill$\Box$

\vskip4mm A ring $R$ is called nil-semicommutative  if $ab=0$ in
$R$ implies that $aRb=0$ for every $a, b \in N(R)$ (see
\cite{MMZ}). For instance, every semicommutative ring (i.e., $ab=0$
in $R$ implies that $aRb=0$) is nil-semicommutative.

\vskip4mm \hspace{-1.8em} {\bf Corollary 3.5.}\ \ {\it Every
nil-semicommutative weakly periodic ring is strongly periodic.}
\vskip2mm\hspace{-1.8em} {\it Proof.} We have from ~\cite[Lemma
2.7]{MMZ} that every nil-semicommutative ring is 2-primal, so the
result follows from Theorem 3.3.\hfill$\Box$

\vskip4mm We note that strongly periodic rings may not be nil-semicommutative as the following shows.

\vskip4mm \hspace{-1.8em} {\bf Example 3.6.}\ \ Let ${\Bbb Z}_2$
be the field of integral modulo $2$, and let
$$R_n=\{
\left(
\begin{array}{ccccc}
a&a_{12}&a_{13}&\cdots&a_{1n}\\
0&a&a_{23}&\cdots&a_{1n}\\
0&0&a&\cdots&a_{3n}\\
\vdots&\vdots&\vdots&\ddots&\vdots\\
0&0&0&\cdots&a
\end{array}
\right)~|~a,a_{ij}\in {\Bbb Z}_2\}$$ with $3\leq n\in {\Bbb N}$.
Let $R=\big( \bigoplus\limits_{n=3}^{\infty}R_n,1\big)$ be the
subalgebra of $\prod\limits_{n=3}^{\infty}R_n$ over ${\Bbb Z}_2$
generated by $\bigoplus\limits_{n=3}^{\infty}R_n$ and $1$. We note
that $P(R)=\bigoplus\limits_{n=3}^{\infty}P(R_n)$. Hence,
$R/P(R)\cong \big( \bigoplus\limits_{n=3}^{\infty}F_n,1\big)$, the
subalgebra of $\prod\limits_{n=3}^{\infty}F_n$ over ${\Bbb Z}_2$
generated by $\bigoplus\limits_{n=3}^{\infty}F_n$ and
$1_{\prod\limits_{n=3}^{\infty}F_n}$, where $F_n={\Bbb Z}_2$ for
all $n=3,4,\cdots $. This implies that $R/P(R)$ is reduced. For
any $a\in N(R)$, $\overline{a}\in R/P(R)$ is nilpotent, and so
$\overline{a}=\overline{0}$. That is, $a\in P(R)$. Therefore $R$
is 2-primal. As $R_n$ is a finite ring for each $n$, we see that
it is periodic. We infer that $R$ is periodic, and so it is weakly
periodic. In light of Theorem 3.3, $R$ is strongly periodic. We
claim that $R_4$ is not nil-semicommutative. Choose
$$a=\left(\begin{array}{ccccc}
0&1&-1&0\\
0&0&0&0\\
0&0&0&0\\
0&0&0&0
\end{array}
\right), x=\left(
\begin{array}{ccccc}
0&0&0&0\\
0&0&1&0\\
0&0&0&0\\
0&0&0&0
\end{array}
\right)~\mbox{and}~b=\left(
\begin{array}{ccccc}
0&0&0&0\\
0&0&0&1\\
0&0&0&1\\
0&0&0&0
\end{array}
\right).$$ Then $a^2=b^2=0$, and so $a,b\in N(R_4)$. Furthermore,
$ab=0$, while $axb\neq 0.$ Thus, $R_4$ is not nil-semicommutative.
Therefore $R$ is not nil-semicommutative, and we are done.

\vskip4mm \hspace{-1.8em} {\bf Theorem 3.7.}\ \ {\it Let $T$ be
the ring of a Morita context $(A,B,$ $M,N,$ $\psi,\varphi)$. If
$im(\psi)$ and $im(\varphi)$ are nilpotent, then $A$ and $B$ are
strongly periodic if and only if so is
$T$.}\vskip2mm\hspace{-1.8em} {\it Proof.}\ \ Suppose $A$ and $B$
are strongly periodic. Then $A$ and $B$ are 2-primal, by Theorem
3.3. Further, they are periodic. In view of Theorem 2.2, $T$ is
periodic. It suffices to prove that $T$ is 2-primal.

Let $\left(
\begin{array}{cc}
a&n\\
m&b
\end{array}
\right)\in T$ is nilpotent. Then we can find some $c\in im(\psi)$
and $d\in im(\varphi)$ such that $a^k+c=0$ and $b^l+d=0$ for some
$k,l\in {\Bbb N}$. This implies that $a\in N(A)$ and $b\in N(B)$.
As $A$ is 2-primal, $a\in P(A)$. Analogously to the proof in Theorem 3.1, we see that $AaA$ is nilpotent.
Likewise, $BbB$ is nilpotent. Clearly, $$T\left(
\begin{array}{cc}
a&n\\
m&b
\end{array}
\right)T\subseteq \left(
\begin{array}{cc}
AaA+im(\psi)&N\\
M&BbB+im(\varphi)
\end{array}
\right).$$ As the sum of two nilpotent ideal of a ring is nilpotent, we see that
$AaA+im(\psi)$ and $BbB+im(\varphi)$ are nilpotent ideals of $A$ and $B$, respectively.
Similarly to the proof of Theorem 2.2, we see that $\left(
\begin{array}{cc}
AaA+im(\psi)&N\\
M&BbB+im(\varphi)
\end{array}
\right)$ is a nilpotent ideal of $T$. Hence, $T\left(
\begin{array}{cc}
a&n\\
m&b
\end{array}
\right)T$ is nilpotent. As in the proof of Theorem 3.1,
we see that $\left(
\begin{array}{cc}
a&n\\
m&b
\end{array}
\right)\in Q$ for any prime ideal $Q$ of $T$. Hence, $\left(
\begin{array}{cc}
a&n\\
m&b
\end{array}
\right)\in P(T)$. Thus, $T$ is 2-primal, and so $T$ is strongly
periodic, by Theorem 3.3.

Conversely, assume that $T$ is strongly periodic. Then $A$ is
periodic. Let $a\in N(R)$. Then $\left(
\begin{array}{cc}
a&0\\
0&0
\end{array}
\right)\in N(T)$. By virtue of Theorem 3.3, $T$ is 2-primal;
hence, $\left(
\begin{array}{cc}
a&0\\
0&0
\end{array}
\right)\in P(T)$. As in the proof in Theorem 3.1, we see that $T\left(
\begin{array}{cc}
a&0\\
0&0
\end{array}
\right)T$ is nilpotent. Then $AaA$ is nilpotent,
and so $a\in P(R)$. It follows that $A$ is 2-primal. Therefore
$A$ is strongly periodic, by Theorem 3.3. Likewise, $B$ is
strongly periodic, as required.\hfill$\Box$

\vskip4mm \hspace{-1.8em} {\bf Corollary 3.8.}\ \ {\it Let $R$ be
strongly periodic, and let $s\in N(R)\bigcap C(R)$. Then
$M_{(s)}(R)$ is strongly periodic.}\vskip2mm\hspace{-1.8em} {\it
Proof.}\ \ As in the proof of Corollary 2.3, we have
$M_s(R)=(R,R,R,R,\psi,\varphi)$ where $im(\varphi)$ and $im(\psi)$
are nilpotent. This completes the proof, by Theorem
3.7.\hfill$\Box$

\vskip4mm \hspace{-1.8em} {\bf Example 3.9.}\ \ Consider the Morita context $R=
\left(
\begin{array}{cc}
{\Bbb Z}_4&{\Bbb Z}_4\\
2{\Bbb Z}_4&{\Bbb Z}_4
\end{array}
\right)$, where the context products are the same as the product in ${\Bbb Z}_4$. Then we claim that $R$ is strongly periodic. Since $R$ is finite, it is periodic, and then we are done by Theorem 3.7.

\vskip4mm As a consequence, a ring $R$ is strongly periodic if and
only if so is the trivial Morita context $M_{(0)}(R)$. Now we exhibit the useful
characterizations of strongly periodic rings as follows.

\vskip4mm \hspace{-1.8em} {\bf Theorem 3.10.}\ \ {\it Let $R$ be a
ring. Then the following are equivalent:}
\begin{enumerate}
\item [(1)] {\it $R$ is strongly periodic.}
\vspace{-.5mm}
\item [(2)] {\it $R/P(R)$ is potent.}
\vspace{-.5mm}
\item [(3)] {\it For any $a\in R$, there exists a prime $m\geq 2$ such that $a-a^m\in P(R)$.}
\vspace{-.5mm}
\item [(4)] {\it For any $a\in R$, $a=eu+w$, where $e=e^2\in R, u^m=1~ (m\in {\Bbb N}),
w\in P(R)$ and $e,u,w$ commutate.}
\end{enumerate}
\vspace{-.5mm} {\it Proof.}\ \ $(1)\Rightarrow (2)$ This is
obvious.

$(2)\Rightarrow (3)$ Luh's Theorem states that a ring $S$ is potent if and only if
for any $x\in S$ there exists a prime $n$ such that $x=x^n$ (cf.~\cite{Lu}). Let $a\in R$. Since $R/J(R)$ is potent,
we have a prime $m\geq 2$ such that
$\overline{a}=\overline{a^m}$ in $R/P(R)$. Therefore, $a-a^m\in
P(R)$.

$(3)\Rightarrow (4)$ Let $a\in R$. Then we have a prime $n\geq 2$
such that $a-a^n\in P(R)\subseteq N(R)$. By Theorem 1.1, $R$ is
periodic. Let $x\in N(R)$. Then $\overline{x}\in R/P(R)$ is
potent; whence, $\overline{x}=\overline{0}$ in $R/P(R)$. Thus,
$x\in P(R)$, and so $N(R)\subseteq P(R)$. By ~\cite[Proposition
13.1.18]{CH1}, $a=eu+w$, where $e=e^2\in R, u\in U(R), w\in P(R)$
and $e,u,w$ commutate. Write $u^{k}=u^{k+m}$ for some $m,k\in
{\Bbb N}$. Then $u^m=1$, as desired,

$(4)\Rightarrow (1)$ For any $a\in R$, $a=eu+w$, where $e=e^2\in
R, u^m=1~ (m\in {\Bbb N}), w\in P(R)$ and $e,u,w$ commutate. Set
$p=eu$. Then $p=eu^{m+1}=p^{m+1}$, i.e., $p\in R$ is potent. Thus,
$R$ is strongly periodic.\hfill$\Box$

\vskip4mm \hspace{-1.8em} {\bf Corollary 3.11.}\ \ {\it Every
subring of a strongly periodic ring is strongly periodic.}
\vskip2mm\hspace{-1.8em} {\it Proof.}\ \ Let $R$ be strongly
periodic, and let $S\subseteq R$. For any $a\in S$, there exists
some $n\geq 2$ such that $a-a^n\in P(R)$ in terms of Theorem 3.10.
Given $a-a^n=a_0,a_1,a_2,\cdots \in S$ with each $a_{i+1}\in
a_iSa_i$, we see that $a-a^n=a_0,a_1,a_2,\cdots \in R$ with each
$a_{i+1}\in a_iRa_i$. This forces that $a_m=0$ for some $m\geq 2$.
Therefore $a-a^n\in P(S)$. By using Theorem 3.10 again, $S$ is
strongly periodic, as needed.\hfill$\Box$

\vskip4mm For example, if $R$ is the finite subdirect product of
strongly periodic rings, then Corollary 3.11 shows that $R$ is
strongly periodic.

\vskip4mm \hspace{-1.8em} {\bf Example 3.12.}\ \ Let $F=GF(q)$ be
a Galois field and let $V$ be an infinite dimensional left vector
space over $F_p$ with $\{v_1,v_2,\cdots \}$ a basis. For the
endomorphism ring $A=End_{F}(V)$, define $A_1=\{ f\in
A~|~rank(f)<\infty ~\mbox{and}~f(v_i)=a_1v_1+\cdots
+a_iv_i~\mbox{for}~i=1,2,\cdots ~\mbox{with}~a_j\in F_p\}$ and let
$R$ be the $F$-algebra of $A$ generated by $A_1$ and $1_A$. Then
$R$ is strongly periodic. By the argument in ~\cite[Example
1.1]{Lee}, $R/P(R)\cong \{ (a_1,\cdots ,a_n,b,b,\cdots
)~|~a_i,b\in F~\mbox{and}~n=1,2,\cdots \}$. As $F=GF(q)$, we see
that $x=x^q$ for all $x\in F$, and then $R/P(R)$ is potent.
According to Theorem 3.10, $R$ is strongly periodic, and we are
through.

\vskip4mm \hspace{-1.8em} {\bf Lemma 3.13.}\ \ {\it Let $I$ be a
nilpotent ideal of a ring $R$. If $R/I$ is strongly periodic, then
so is $R$.} \vskip2mm\hspace{-1.8em} {\it Proof.}\ \ Let $a\in R$.
Then there exists some $n\geq 2$ such that $\overline{a-a^n}\in
P\big(R/I\big)$. Hence, $\big(R(a-a^n)R\big)^m\subseteq I$. As $I$
is nilpotent, $\big(R(a-a^n)R\big)^{mn}=0$. This shows that
$a-a^n\in P(R)$. Therefore $R$ is strongly periodic, by Theorem
3.9.\hfill$\Box$

\vskip4mm \hspace{-1.8em} {\bf Theorem 3.14.}\ \ {\it Let $I$ be
an ideal of a ring $R$. Then the following are equivalent:}
\begin{enumerate}
\item [(1)] {\it $R/I$ is strongly periodic.}
\vspace{-.5mm}
\item [(2)] {\it $R/I^n$ is strongly periodic for all $n\in {\Bbb N}$.}
\vspace{-.5mm}
\item [(2)] {\it $R/I^n$ is strongly periodic for some $n\in {\Bbb N}$.}
\end{enumerate}
\vspace{-.5mm} {\it Proof.}\ \ $(1)\Rightarrow (2)$ Clearly,
$R/I\cong \big(R/I^n\big)/\big(I/I^n\big)$. Since
$\big(I/I^n\big)^n=0$, proving $(2)$ by Lemma 3.13.

$(2)\Rightarrow (3)$ This is trivial.

$(3)\Rightarrow (1)$ For any $\overline{a}\in R/I$, we see that
$a+I^n\in R/I^n$. By hypothesis, there exists a potent
$\overline{p}\in R/I^n$ such that $\overline{a-p}\in
P\big(R/I^n\big)$. Write $\overline{p}=\overline{p}^m$ for some
$m\geq 2$. Then $p-p^m\in I^n\subseteq I$, and so $\overline{p}\in
R/I$ is potent. Obviously, $(R/I^n)\overline{(a-p)}(R/I^n)$ is
nilpotent, and then $\big(R(a-p)R\big)^s\subseteq I^n\subseteq I$
for some $s\in {\Bbb N}$. We infer that
$(R/I)\overline{(a-p)}(R/I)$ is nilpotent. As in the proof of Theorem 3.1, we infer that
$\overline{a-p}\in P(R/I)$, as required.\hfill$\Box$

\vskip4mm Recall that a ring $R$ is an abelian
ring if every idempotent in $R$ is central. A ring $R$ is strongly $\pi$-regular if for any $a\in R$ there exists $n\in {\Bbb N}$ such that
$a^n\in a^{n+1}R$. Obviously, every periodic ring is strongly $\pi$-regular. We now derive

\vskip4mm \hspace{-1.8em} {\bf Lemma 3.15.}\ \ {\it Every abelian
periodic ring is strongly periodic.} \vskip2mm\hspace{-1.8em} {\it
Proof.}\ \ Let $R$ be an abelian periodic ring. Then $R$ is
strongly $\pi$-regular. Badawi's Theorem states that the set of all nilpotent elements of an abelain strongly $\pi$-regular ring is an ideal (~\cite{BA}). Thus, $N(R)$ forms an ideal of $R$. This completes the proof, by Theorem 3.1.\hfill$\Box$

\vskip4mm Let $n\geq 2$ be a fixed integer. A ring $R$ is said to
be generalized $n$-like provided that for any $a,b\in R$,
$(ab)^n-ab^n-a^nb+ab=0$ (cf.~\cite{M}). It is proved that every generalized $3$-like ring is commutative (~\cite[Theorem 3]{M}). We now derive

\vskip4mm \hspace{-1.8em} {\bf Theorem 3.16.}\ \ {\it Every
generalized $n$-like ring is strongly periodic.}
\vskip2mm\hspace{-2.0em} {\it Proof.}\ \ Let $R$ be a generalized $n$-like ring, and let $a\in R$. Then
$a^{2n}-2a^{n+1}+a^2=0$, and so $(a-a^n)^2=0$. Thus, $a-a^n\in
N(R)$. Accordingly, $R$ is periodic by Theorem 1.1. In light of ~\cite[Lemma 2]{M}, $R$ is
abelian. Therefore $R$ is strongly periodic, by Lemma 3.15.\hfill$\Box$

\vskip4mm Let $R=\{
\left( \begin{array}{ccc} x&y&z\\
0&x^2&0\\
0&0&x
\end{array}
\right)~|~x,y,z\in GF(4)\}$. Then for each $a\in R$, $a^7=a$ or
$a^7=a^2=0$. Therefore $R$ is a generalized $7$-like ring. By
Theorem 3.16, $R$ is strongly periodic. In this case, $R$ is abelian but not commutative (cf.~\cite[Example 2]{M}).

\section{J-Clean-like Rings}

\vskip4mm We now consider J-clean-like Morita contexts and extend
Theorem 2.2 as follows.

\vskip4mm \hspace{-1.8em} {\bf Theorem 4.1.}\ \ {\it Let $T$ be
the ring of a Morita context $(A,$ $B,$ $M,$ $ N,$ $
\psi,\varphi)$ with $im(\psi)\subseteq J(A)$ and
$im(\varphi)\subseteq J(B)$. If $A$ and $B$ are J-clean-like, then
so is $T$.}\vskip2mm\hspace{-1.8em} {\it Proof.}\ \ Let $\left(
\begin{array}{cc}
a&n\\
m&b
\end{array}
\right)\in T.$ Then we have potent $p\in A$ and $q\in B$ such that
$a-p\in J(A)$ and $b-q\in J(B$.
Hence
$$
\left(
\begin{array}{cc}
a&n\\
m&b
\end{array}
\right) - \left(
\begin{array}{cc}
p&0\\
0&q
\end{array}
\right) = \left(
\begin{array}{cc}
a-p&n\\
m&b-q
\end{array}
\right).$$ Let $\left(
\begin{array}{cc}
c&s\\
t&d
\end{array}
\right)\in T$. As $1_A-(a-p)c-\psi (n\bigotimes t)\in U(A)$ and $1_B-(b-q)d-\varphi(m\bigotimes s)\in U(B)$, it follows by ~\cite[Lemma 3.1]{TZ} that $$\begin{array}{lll}
&&1_T-\left(
\begin{array}{cc}
a-p&n\\
m&b-q
\end{array}
\right)\left(
\begin{array}{cc}
c&s\\
t&d
\end{array}
\right)\\
&=&\left(
\begin{array}{cc}
1_A-(a-p)c-\psi (n\bigotimes t)&*\\
**&1_B-(b-q)d-\varphi(m\bigotimes s)
\end{array}
\right)\in U(T).\end{array}$$ Hence, $\left(
\begin{array}{cc}
a-p&n\\
m&b-q
\end{array}
\right)\in J(T)$, and therefore $T$ is J-clean-like.\hfill$\Box$

\vskip4mm As a consequence, we deduce that the $n\times n$ lower
(upper) triangular matrix ring over a J-clean-like ring is
J-clean-like.

\vskip4mm \hspace{-1.8em} {\bf Corollary 4.2.}\ \ {\it Let $R$ be
J-clean-like, and let $s\in J(R)\bigcap C(R)$. Then $M_{(s)}(R)$
is J-clean-like.} \vskip2mm\hspace{-1.8em} {\it Proof.}\ \ As in
the proof of Corollary 2.3, $M_{(s)}(R)$ can be regarded as the
ring of a Morita context $(R,R,R,R,\psi,\varphi)$ with
$im(\psi)\subseteq J(R)$ and $im(\varphi)\subseteq J(R)$.
According to Theorem 4.1, $M_{(s)}(R)$ is
$J$-clean-like.\hfill$\Box$

\vskip4mm \hspace{-1.8em} {\bf Corollary 4.3.}\ \ {\it Let $R$ be
a J-clean-like ring. Then $M_{(x)}(R[[x]])$ is J-clean-like.}
\vskip2mm\hspace{-1.8em} {\it Proof.}\ \ For any $f(x)\in R[[x]]$,
we can find an potent $p\in R$ such that $f(0)-p\in J(R)$. Hence,
$f(x)=p+\big(f(x)-p\big)$. One easily checks that $f(x)-p\in
J\big(R[[x]]\big)$. Thus, $R[[x]]$ is J-clean-like. Choose $s=x$.
Applying Corollary 4.2 to $R[[x]]$, $M_{(x)}(R[[x]])$ is
$J$-clean-like.\hfill$\Box$

\vskip4mm Analogously, if $R$ is a J-clean-like ring then so is
$M_{(x^m)}\big(R[[x]]/(x^n)\big)$ for all $1\leq m\leq n$.

\vskip4mm \hspace{-1.8em} {\bf Proposition 4.4.}\ \ {\it A ring
$R$ is strongly periodic if and only if the following two conditions hold simultaneously:}
\begin{enumerate}
\item [(1)] {\it $R$ is J-clean-like;}
\vspace{-.5mm}
\item [(2)] {\it $J(R)$ is locally nilpotent.}
\end{enumerate}
\vspace{-.5mm} {\it Proof.}\ \ $\Longrightarrow $ Suppose $R$ is
strongly periodic. As $P(R)\subseteq J(R)$, $R$ is J-clean-like.
Let $x\in J(R)$. Then there exists a potent $p\in R$ such that
$x-p\in P(R)$; hence, $p=x-(x-p)\in J(R)$. This shows that $p=0$,
and so $x\in P(R)$. As in the proof of Theorem 3.1, $RxR$ is nilpotent. As the sum of finite nilpotent ideal is nilpotent, we prove that
$J(R)$ is locally nilpotent, as required.

$\Longleftarrow $ Let $x\in J(R)$. Since $J(R)$ is locally nilpotent, $RxR$ is nilpotent. As in the proof of Theorem 3.1, we get $x\in P(R)$. Hence, $J(R)\subseteq P(R)$. This completes the proof, by
$(1)$.\hfill$\Box$

\vskip4mm Recall that a ring $R$ is J-clean provided the for any
$a\in R$ there exists an idempotent $e\in R$ such that $a-e\in
J(R)$ (cf. ~\cite{CH}). This following result explains the
relation between J-clean rings and J-clean-like rings.

\vskip4mm \hspace{-1.8em} {\bf Proposition 4.5.}\ \ {\it A ring
$R$ is J-clean if and only if the following two conditions hold:}
\begin{enumerate}
\item [(1)] {\it $R$ is J-clean-like;}
\vspace{-.5mm}
\item [(2)] {\it $J(R)=\{ x\in R~|~1-x\in U(R)$\}.}
\end{enumerate}
\vspace{-.5mm} {\it Proof.}\ \ $\Longrightarrow $ Clearly, $R$ is
J-clean-like. It is easy to check that $J(R)\subseteq \{ x\in
R~|~1-x\in U(R)\}$. If $1-x\in U(R)$, then there exists an
idempotent $e\in R$ such that $w:=x-e\in J(R)$. Hence,
$1-e=(1-x)+w=(1-x)\big(1+(1-x)^{-1}w\big)\in U(R)$. This show that
$1-e=1$, and so $e=0$. Therefore $x\in J(R)$, and so
$J(R)\supseteq \{ x\in R~|~1-x\in U(R)\}$, as required.

$\Longleftarrow$ For any $a\in R$ there exists a potent $p\in R$
such that $(a-1)-p\in J(R)$. Write $p=p^m (m\geq 2)$. Then
$p^{m-1}\in R$ is an idempotent. Set $e=1-p^{m-1}$ and
$u=p-1+p^{m-1}$. Then $e=e^2\in R$ and
$u^{-1}=p^{m-1}-1+p^{m-1}p^{m-2}$. Further, $p=e+u$. This shows
that $a-1=p+(a-1-p)=e+u+(a-1-p)$. Hence, $a=e+\big(u+(a-p)\big)$. As
$1-\big(u+(a-p)\big)=-u-(a-1-p)=-u\big(1-u^{-1}(a-1-p)\big)\in
U(R)$, we see that $u+(a-p)\in J(R)$. Therefore $R$ is J-clean,
as asserted.\hfill$\Box$

\vskip4mm \hspace{-1.8em} {\bf Example 4.6.}\ \ Let
$R=\left( \begin{array}{cc} {\Bbb Z}_3&{\Bbb Z}_3\\
0&{\Bbb Z}_3 \end{array} \right)$. Then $R$ is J-clean-like, while
it is not J-clean.
For any $\left( \begin{array}{cc} a&c\\
0&b \end{array} \right)\in R$, we see that $\left( \begin{array}{cc} a&c\\
0&b \end{array} \right)=\left( \begin{array}{cc} a&0\\
0&b \end{array} \right)+\left( \begin{array}{cc} 0&c\\
0&0 \end{array} \right)$ is the sum of a potent element in $R$ and
an element in $J(R)$, hence that $R$ is J-clean-like. As
$R/J(R)\cong {\Bbb Z}_3$ is not Boolean, we conclude that $R$ is
not J-clean.

\vskip4mm An element $p\in R$ is $J$-potent provided that there
exists some $n\geq 2$ such that $p-p^n\in J(R)$. We say that every
potent element lifts modulo $J(R)$ if for any $J$-potent $p\in R$
there exists a potent $q\in R$ such that $p-q\in J(R)$.

\vskip4mm \hspace{-1.8em} {\bf Lemma 4.7.}\ \ {\it A ring $R$ is
J-clean-like if and only if the following two conditions hold:}
\begin{enumerate}
\item [(1)] {\it $R/J(R)$ is potent;}
\vspace{-.5mm}
\item [(2)] {\it Every potent element lifts modulo $J(R)$.}
\end{enumerate}
\vspace{-.5mm} {\it Proof.}\ \  $\Longrightarrow $ This is
obvious.

$\Longleftarrow $ Let $a\in R$. Then $\overline{a}\in R/J(R)$ is
potent. By hypothesis, we can find a potent $p\in R$ such that
$a-p\in J(R)$. Accordingly, $R$ is J-clean-like.\hfill$\Box$

\vskip4mm Recall that a ring $R$ is right (left) quasi-duo
provided that every maximal right (left) ideal is a two-sided
ideal. As is well known, every right (left) duo ring (i.e., every
right (left) ideal is two-sided) is right (left) quasi-duo. We
come now to the main result of this section.

\vskip4mm \hspace{-1.8em} {\bf Theorem 4.8.}\ \ {\it A ring $R$ is
J-clean-like if and only if the following three conditions hold:}
\begin{enumerate}
\item [(1)] {\it $R/J(R)$ is periodic;}
\vspace{-.5mm} \vspace{-.5mm}
\item [(2)] {\it $R$ is right (left) quasi-duo;}
\item [(3)] {\it Every potent element lifts modulo $J(R)$.}
\end{enumerate}
\vspace{-.5mm} {\it Proof.}\ \ $\Longrightarrow$ In view of Lemma
4.7, $R/J(R)$ is potent, and so it is periodic. Let $M$ be a
maximal right ideal of $R$, and let $r\in R$. Then $J(R)\subseteq
M$, and that $M/J(R)$ is a maximal right of $R/J(R)$. As is well
known, every potent ring is commutative, and so $rx+J(R)\in
M/J(R)$ for any $x\in M$. Write $rx+J(R)=y+J(R)$ for a $y\in M$.
hence, $rx-y\in J(R)\subseteq M$. This shows that $rx\in M$;
hence, $rM\subseteq M$. Therefore $M$ is a two-sided ideal, and
then $R$ is right quasi-duo. Likewise, $R$ is left quasi-duo.
$(3)$ is obvious, by Lemma 4.7.

$\Longleftarrow $ Let $a\in R$. By $(1)$, there exists a $p\in R$ such
that $\overline{a-p}\in N\big(R/J(R)\big), p-p^n\in J(R) (n\geq
2)$, by Theorem 1.1. By $(3)$, we may assume that $p=p^n$. Set
$w=a-p$. Then $\overline{w}^m=0$. Since $R$ is right (left)
quasi-duo, as in ~\cite[Corollary 3.4.7]{CH1}, we see that
$ex-xe\in J(R)$ for any idempotent $e\in R$ and any element $x\in R$.
This means that $R/J(R)$ is
abelian. Similarly to the proof of ~\cite[Corollary 1.3.15]{CH1},
$R/J(R)$ is reduced. Hence, $\overline{w}=\overline{0}$, and then
$w\in J(R)$. Therefore $a-p\in J(R)$, as desired.\hfill$\Box$

\vskip4mm As an consequences of Corollary 3.5, every right (left)
duo periodic ring is strongly periodic. Further, we derive

\vskip4mm \hspace{-1.8em} {\bf Corollary 4.9.}\ \ {\it A ring $R$
is strongly periodic if and only if }
\begin{enumerate}
\item [(1)] {\it $R$ is periodic;}
\vspace{-.5mm} \vspace{-.5mm}
\item [(2)] {\it $R$ is right (left) quasi-duo;}
\item [(3)] {\it $J(R)$ is locally nilpotent.}
\end{enumerate}
\vspace{-.5mm} {\it Proof.}\ \  $\Longrightarrow$ Clearly, $R$ is
periodic. It follows from Proposition 4.4 that $R$ is J-clean-like
and $J(R)$ is locally nilpotent. Thus, $R$ is right (left)
quasi-duo, by Theorem 4.8.

$\Longleftarrow $ Since $R$ is periodic, $R/J(R)$ is periodic.
Thus, $R$ is J-clean-like, by Theorem 4.8. By $(3)$, $J(R)=P(R)$,
and the result follows.\hfill$\Box$

\vskip4mm \hspace{-1.8em} {\bf Example 4.10.}\ \ {\it Let $R={\Bbb
Z}_{(5)}$. Then $R$ is right (left) quasi-duo, $R/J(R)$ is
periodic, while $R$ is not $J$-clean-like.}
\vskip2mm\hspace{-1.8em} {\it Proof.}\ \ Let $R={\Bbb Z}_{(5)}$.
Then $J(R)=5R$. Hence, $R/J(R)\cong {\Bbb Z}_{5}$ is a finite
field. Thus, $R/J(R)$ is periodic. Suppose every potent element
lifts modulo $J(R)$. Clearly, $2-2^5\in J(R)$. Hence,
$\overline{2}\in R/J(R)$ is potent. Thus, we can find a potent
$w\in R$ such that $2-w\in J(R)$. Write $w=\frac{m}{n}$, where
$(m,n)=1,5\nmid n$ and $w=w^s(s\geq 2)$. Then $w(1-w^{s-1})=0$,
and so $w=0$ or $w^{s-1}=1$. If $w=0$, then $2\in J(R)$, a
contradiction. If $w^{s-1}=1$, then $\frac{m^{s-1}}{n^{s-1}}=1$;
whence, $m=\pm n$. This implies that $w=\pm 1$; hence,
$2-w=1,3\not\in J(R)$, a contradiction. Therefore $R$ is not
J-clean-like , by Lemma 4.7.\hfill$\Box$

\vskip4mm \hspace{-1.8em} {\bf Lemma 4.11.}\ \ {\it Let $R$ be
J-clean-like. Then $N(R)\subseteq J(R)$.} \vskip2mm\hspace{-1.8em}
{\it Proof.}\ \ Let $x\in N(R)$. Then $x^m=0$ for some $m\geq 2$.
Moreover, there exists a potent $p\in R$ such that $w:=x-p\in
J(R)$. Write $p=p^n$ for some $n\geq 2$. Then
$p=p^n=(p^n)^n=p^{n^2}=(p^n)^{n^2}=p^{n^3}=\cdots =p^{n^m}$.
Clearly, $n^m=\big(1+(n-1)\big)^m\geq m(n-1)\geq m$, and so
$x^{n^m}=0$. As $x^{n^m}-p^{n^m}\in J(R)$, we have
$x=p+w=p^{n^m}+w\in J(R)$. Therefore $x\in J(R)$, hence the
result.\hfill$\Box$

\vskip4mm \hspace{-1.8em} {\bf Lemma 4.12.}\ \ {\it Let $R$ be a
ring. Then the following are equivalent:}
\begin{enumerate}
\vspace{-.5mm} \item [(1)] {\it $R$ is a periodic ring in which
every nilpotent is contained in $J(R)$.} \vspace{-.5mm}
\item [(2)] {\it $R$ is J-clean-like and $J(R)$ is nil.}
\end{enumerate}
\vspace{-.5mm}  {\it Proof.}\ \ $(1)\Rightarrow (2)$ Suppose $R$
is a periodic ring with $N(R)\subseteq J(R)$. Let $x\in J(R)$. Then we have $m,n\in {\Bbb N}$ such that $x^m=x^n (n>m)$. Hence,
$x^m(1-x^{n-m})=0$, and so $x^m=0$. This shows that $J(R)$ is nil. Let $a\in R$. In view of
Theorem 1.1, there exists a potent $p\in R$ such that $a-p\in
N(R)$. By hypothesis, $a-p\in J(R)$. Therefore $R$ is
J-clean-like.

$(2)\Rightarrow (1)$ For any $a\in R$, there exists a potent $p\in
R$ such that $w:=a-p\in J(R)$. Hence, $a=p+w$ and $p=p^n$ for some
$n\geq 2$. Thus, $a^n=p^n+v$ for a $v\in J(R)$. This implies that
$a-a^n=w-v\in J(R)\subseteq N(R)$. Therefore $R$ is periodic, by
Theorem 1.1. In light of Lemma 4.11, every nilpotent of $R$ is
contained in $J(R)$, as desired.\hfill$\Box$

\vskip4mm \hspace{-1.8em} {\bf Theorem 4.13.}\ \ {\it Let $R$ be a
ring. If for any sequence of elements $\{ a_i\}\subseteq R$ there
exists a $k\in {\Bbb N}$ and $n_1,\cdots ,n_k\geq 2$ such that
$(a_1-a_1^{n_1})\cdots (a_k-a_k^{n_k})=0$, then $R$ is
$J$-clean-like.} \vskip2mm\hspace{-1.8em} {\it Proof.}\ \ For any
$a\in R$, we have a $k\in {\Bbb N}$ and $n_1,\cdots ,n_k\geq 2$
such that $(a-a^{n_1})\cdots (a-a^{n_k})=0$. This implies that
$a^k=a^{k+1}f(a)$ for some $f(t)\in {\Bbb Z}[t]$. In view of
Theorem 1.1, $R$ is periodic.

Clearly, $R/J(R)$ is isomorphic to a subdirect product of some
primitive rings $R_i$. Case 1. There exists a subring $S_i$ of
$R_i$ which admits an epimorphism $\phi_i: S_i\to M_{2}(D_i)$
where $D_i$ is a division ring. Case 2. $R_i\cong M_{k_i}(D_i)$
for a division ring $D_i$. Clearly, the hypothesis is inherited by
all subrings, all homomorphic images of $R$, we
claim that, for any sequence of elements $\{ a_i\}\subseteq
M_2(D_i)$ there exists $s\in {\Bbb N}$ and $m_1,\cdots ,m_s\geq 2$
such that $(a_1-a_1^{m_1})\cdots (a_s-a_s^{m_s})=0$. Choose
$a_i=e_{12}$ if $i$ is odd and $a_i=e_{21}$ if $i$ is even. Then
$(a_1-a_1^{m_1})(a_2-a_2^{m_2})\cdots (a_s-a_s^{m_s})=a_1a_2\cdots
a_s\neq 0$, a contradiction. Thus, Case I do not happen. Further, in Case II, $k_i=1$ for all $i$. This shows that
each $R_i$ is reduced, and then so is $R/J(R)$. If $a\in
N(R)$, we have some $n\in {\Bbb N}$ such that $a^n=0$, and thus
$\overline{a}^n=0$ is $R/J(R)$. Hence, $\overline{a}\in
J\big(R/J(R)\big)=0$. This implies that $a\in J(R)$, and so
$N(R)\subseteq J(R)$. Therefore $R$ is J-clean-like , by Lemma
4.12.\hfill$\Box$

\vskip4mm Recall that a subset $I$ of a ring $R$ is left (resp.,
right) $T$-nilpotent in case for every sequence $a_1,a_2,\cdots $
in $I$ there is an $n$ such that $a_1\cdots a_n=0$ (resp.,
$a_n\cdots a_1=0$). Every nilpotent ideal is left and right
$T$-nilpotent. The Jacobson radical $J(R)$ of a ring $R$ is left
(resp., right) $T$-nilpotent if and only if for any nonzero left
(resp., right) $R$-module $M$, $J(R)M\neq M$(resp., $MJ(R)\neq
M$).

\vskip4mm \hspace{-1.8em} {\bf Corollary 4.14.}\ \ {\it Let $R$ be
a ring. If $R/J(R)$ is potent and $J(R)$ is left (resp., right)
$T$-nilpotent, then $R$ is J-clean-like.} \vskip2mm\hspace{-1.8em}
{\it Proof.}\ \ We may assume $R/J(R)$ is potent and $J(R)$ is
left $T$-nilpotent. For every sequence $a_1,a_2,\cdots, a_m,\cdots
$ in $R$, there exists some $n_i\in {\Bbb N}$ such that
$a_i-a_i^{n_i}\in J(R)$ for all $i$. We choose $b_1=a_1-a_1^{n_1},
b_2=(1-b_1)^{-1}\big(a_2-a_2^{n_2}\big),
b_3=(1-b_2)^{-1}\big(a_3-a_3^{n_3}\big),\cdots
,b_m=(1-b_{m-1})^{-1}\big(a_m-a_m^{n_m}\big),\cdots .$ By
hypothesis, we can find some $k\in {\Bbb N}$ such that
$b_1(1-b_1)b_2(1-b_2)\cdots b_{k-1}(1-b_{k-1})=0$. Hence,
$b_1(1-b_1)b_2(1-b_2)\cdots b_{k-1}(1-b_{k-1})b_k=0$. This shows
that $(a_1-a_1^{n_1})\cdots (a_s-a_k^{n_k})=0$. Therefore $R$ is
J-clean-like , by Theorem 4.13.\hfill$\Box$

\vskip6mm \bc{\bf Acknowledgements}\ec

\vskip4mm The authors are grateful to the referee for his/her
helpful suggestions which make the new version clearer.

\end{document}